\documentclass[12pt]{amsart}
\usepackage{mathpazo,fullpage,relsize}
\usepackage{amssymb}
\usepackage{hyperref}
\usepackage{pictexwd,dcpic}
\usepackage{graphicx}
\usepackage{pinlabel}
\usepackage{csquotes}
\newtheorem{thm}{Theorem}
\newtheorem{prop}[thm]{Proposition}
\newtheorem{cor}[thm]{Corollary}
\newtheorem{lem}[thm]{Lemma}       

\newtheorem{conj}[thm]{Conjecture}

\theoremstyle{remark}

\newtheorem{ex}[thm]{Example}

\theoremstyle{definition}
\newtheorem{defn}[thm]{Definition}

\newcommand{\R}{\mathbb{R}}
\newcommand{\Z}{\mathbb{Z}}
\newcommand{\C}{\mathbb{C}}

\title[{Doubles without open book decompositions from higher signatures}]{Doubles without open book decompositions\\ from higher signatures}
\author{D.~Kotschick}
\address{Mathematisches Institut, {\smaller LMU} M\"unchen, Theresienstr.~39, 80333~M\"unchen, Germany}
\email{dieter@math.lmu.de}
\date{\today; \copyright{\ D.~Kotschick 2025}}
\thanks{This material is based upon work supported by a grant from the Institute
for Advanced Study School of Mathematics in Princeton.}

\begin{document}

\begin{abstract}
We show that in every even dimension there are closed manifolds that are doubles, but 
have no open book decomposition. In high dimensions, this contradicts the conclusions in Ranicki's 
book~\cite{Ran}. In all dimensions, examples arise from the non-multiplicativity of the signature in 
fibre bundles. We discuss many examples and applications in dimension four, where this 
phenomenon is related to the simplicial volume.
\end{abstract}

\maketitle


\section{Introduction}

The concept of an open book was introduced by Winkelnkemper~\cite{Win} in 1973.
Around the same time, Tamura~\cite{Tam} introduced spinnable structures, which are 
in fact equivalent to open books, but it is Winkelnkemper's terminology that has become 
standard. These authors noted that  a result of Alexander from 1923 can be restated
as saying that every three-manifold has the structure of an open book. The work
of Winkelnkemper~\cite{Win}, Tamura~\cite{Tam}, T.~Lawson~\cite{TLaw} and Quinn~\cite{Quinn}
extended this to all odd dimensional manifolds.

The even dimensional case is more complicated. Winkelnkemper~\cite{Win} already 
noted that the signature is an obstruction to the existence of an open book, and he 
stated an existence theorem for open book decompositions of simply-connected 
manifolds with vanishing signature in dimensions $\geq 8$. This was generalised
to arbitrary fundamental groups, and extended to dimension $6$, by Quinn~\cite{Quinn},
who proved that the existence of an open book decomposition of $M$ is equivalent to the vanishing of a 
certain invariant $i(M)\in W^s(\Z[\pi_1(M)])$ with values in a Witt group of sesquilinear
forms over the group ring of the fundamental group. For trivial fundamental groups,
the Quinn invariant reduces to the signat                                                                                                                                                                                                                                                                                                                                                                                                                                                                                                                                                            ure.

Winkelnkemper was well aware of a relationship between doubles and open book
decompositions; cf.~Section~\ref{s:ob} below. 
Quinn also mentioned this, in fact his paper~\cite{Quinn} contains the following sentence at the top of page~57:
                                                                                                                                                                                                                                                                                                                                                                                                                                                                                                                                                                                                                                                           
\begin{displayquote}
{\it Book decompositions are closely related to twisted doubles (book $\Rightarrow$ double, and perhaps conversely).}
\end{displayquote}

Clearly, in ``book $\Rightarrow$ double'' he meant twisted doubles. The main point of this paper
is to observe that the converse of this implication, which Quinn hinted at, but was careful 
not to claim, is false. Moreover, the counterexamples are ordinary rather than twisted doubles.
They are detected by the well known non-multiplicativity of the signature in fibre bundles, discovered
by Atiyah~\cite{A} and Kodaira~\cite{Kd}. In view of the result of Chern, Hirzebruch and Serre~\cite{CHS},
showing that the signature is in fact multiplicative if the fundamental group of the base acts trivially
on the cohomology of the fibre, this phenomenon is tied to the fundamental group, and in 
particular to local coefficient systems arising from its representations. The Quinn invariant 
has to capture all the possible signatures of all coefficient systems over the manifold.

There is an account of open book decompositions, respectively twisted doubles, 
in Ranicki's book~\cite{Ran} in Chapter~29, respectively Chapter~30. Unfortunately 
the conclusions reached there are not correct. Ranicki in particular wrote on page~371:

\begin{displayquote}
{\it Thus a high-dimensional manifold is a twisted double if and only if it is an open book.}
\end{displayquote}

So Ranicki wrongly claimed the converse to the implication Quinn had mentioned.
This has led to some confusion in the literature on open book decompositions, to which I was 
alerted by the preprint of Kastenholz~\cite{Kas}.

In this paper I develop a sequence of examples of manifolds that are doubles yet have no 
open book decomposition. They are detected simply by the fact that there are surface bundles 
over surfaces with non-zero signature, and I like to think of this argument as being completely
elementary. It has turned out a posteriori that such a discussion can be phrased in the 
context of the so-called SK groups, as was done by Neumann~\cite{Neu}. I have chosen 
not to do this, in order not to obscure the very simple geometric mechanism that shows that
the conclusions in~\cite{Ran} are wrong. Some further comments on this are contained 
in the final section of the paper.

I also discuss the relationship of these elementary arguments to the notion of simplicial 
volume, as used by Kastenholz~\cite{Kas}. This turns out to be special to dimensions $2$ and $4$.
Conjecturally, in these dimensions only, the non-vanishing of obstructions to open books arising from the 
non-multiplicativity of the signature are in fact equivalent to the positivity of the simplicial 
volume in the sense of Gromov~\cite{Gro}.

In the final section we show that there are many Engel four-manifolds without open book decompositions. 
Since, based on~\cite{Kas}, this conclusion was also reached independently by Lawande and Saha~\cite{LS},
we keep our account brief.

\subsection*{Acknowledgements} 
This paper started out as a letter to T.~Kastenholz. I am grateful to him for his correspondence, and to 
M.~Stover for an exchange about Example~\ref{Stover}.

\section{Open books and twisted doubles}\label{s:ob}

All the manifolds $M$ considered in this paper are smooth and oriented. If we flip the (implicit)
orientation of $M$, the resulting oriented manifold is denoted by $\overline{M}$.

There are several different ways one can define open books. 
The following is essentially Win\-keln\-kem\-per's definition from~\cite{Win}, phrased carefully as
done for example by Tamura~\cite{Tam}, Durfee--Lawson~\cite{DL} or Giroux~\cite{Gir}:

\begin{defn}\label{GOB}
An {\it open book decomposition} of a smooth manifold $M$  consists of a smooth codimension $2$
submanifold $B$ with trivial normal bundle, so that there are tubular neighbourhoods $N$ of $B$ in $M$ of the form
$B\times D^2$, and a locally trivial smooth fibration 
$$
\theta\colon M\setminus B \longrightarrow S^1
$$
such that for a suitable neighbourhood $N$ as above, the restriction of $\theta$ to $N\setminus B$
is given by the composition of obvious projections
$$
\theta\vert_{N\setminus B}\colon N\setminus B = 
B\times (D^2\setminus \{ 0\}) \stackrel{\pi_2}{\longrightarrow} D^2\setminus \{ 0\} \stackrel{z/\vert z\vert}{\longrightarrow} S^1 \ .
$$
\end{defn}

The following alternative definition was given by H.~B.~Lawson~\cite{Law}. The name of Alexander appears because 
of his result about three-manifolds mentioned in the introduction.
\begin{defn}\label{Alex}
A smooth manifold $M$ has an {\it Alexander decomposition} if there is a smooth 
function $p\colon M\longrightarrow\C$ such that
\begin{itemize}
\item[(1)] $0\in\C$ is a regular value of $p$, 
and 
\item[(2)] $\frac{p}{\vert p\vert} \colon M\setminus p^{-1}(0) \longrightarrow S^1$ is a submersion.
\end{itemize}
\end{defn}
It is clear that this is equivalent to an open book decomposition. A detailed proof of the equivalence was given by Schiegnitz~\cite{Sch}.
The first condition implies that $B = p^{-1}(0)$ is a smooth submanifold of codimension $2$, and that 
$p$ is a submersion on a neighbourhood of this submanifold. In particular $B$ has trivial normal bundle in $M$.
The second condition implies that the complement of $B$ fibers smoothly over the circle.

The submanifold $B$ is called the {\it binding} of the open book decomposition, and the fibres of $\theta$ or $p/\vert p\vert$
are the {\it pages}, all diffeomorphic to the interior of a compact manifold $P$ with boundary $\partial P =B$.
The open book decomposition is determined by the monodromy of the fibration of $M\setminus N$ over 
the circle, which is a diffeomorphism $f\colon P\longrightarrow P$ restricting to the identity near $\partial P = B$.
In keeping with our orientation assumption, we consider only the case in which $P$ is oriented. We orient $B$ as the 
boundary of $P$. Since $f$ is the identity on a neighbourhood of the boundary, it is necessarily orientation-preserving,
and the orientation of $M$ can be recovered from the orientation of $P$ and that of $S^1$.

Following Bowden--Crowley~\cite{BC}, a twisted double decomposition of an open book can be described as follows.
The union of two distinct pages, attached to each other at the binding, is a copy of the double $D(P)=P\cup \overline{P}$
embedded in $M$ as a smooth separating hypersurface $M$. The closures of the connected components of the 
complement are copies of $P \times [0,1]$ (up to smoothing corners). This shows that $M$ is a twisted
double obtained by gluing two copies of $P \times [0,1]$ along the boundary $D(P)=P\cup \overline{P}$ by the diffeomorphism 
$f\cup Id_{\overline{P}}$.

If the monodromy diffeomorphism $f$ in an open book decomposition of $M$ can be arranged to be the identity,
then $M$ is an honest double, rather than a twisted double.

\section{Non-existence of open books}\label{s:noOB}

The following three Lemmata are quite standard, see for example~\cite{Sch}.
\begin{lem}\label{Euler}
If $M$ admits an open book decomposition with binding $B$, then the Euler characteristics of $M$ and $B$ coincide.
\end{lem}
\begin{proof}
This is a straightforward application of the Mayer-Vietoris argument.
\end{proof}

\begin{lem}\label{Weyl}
If $M$ admits an open book decomposition, then its signature vanishes.
\end{lem}
\begin{proof}
This is usually described as a consequence of Novikov additivity, but can also be seen by the following completely
different argument, inspired by~\cite{MPCPS,K}.

Consider the collection of all manifolds with open book decompositions with a fixed page $P$ and 
arbitrary monodromy. The Betti numbers of these manifolds are bounded independently of the monodromy
through the Betti numbers of $P$, and so the signature, which is bounded in absolute value by the middle-degree
Betti number, is also bounded. We will see that this forces it to vanish.

A given monodromy diffeomorphism $f$ gives a particular manifold $M_1=M(P,f)$,
and iterating $f$ produces $M_k=M(P,f^k)$ which is a $k$-fold cyclic cover of $M_1$ branched in $B$.
In general, the signature is not always multiplicative in branched covers, since there is a correction term involving 
the branch locus, see Hirzebruch~\cite{Hir}. However, this correction term depends only on the homology class of the branch locus.
Here the branch locus is the binding, and so, being the boundary of the page $P$, is homologically trivial. 
Therefore, the signature is in fact multiplicative: $\sigma (M_k) = k \cdot \sigma (M_1)$. Since this is bounded independently of $k$,
we conclude $\sigma (M_1)=0$.
\end{proof}

\begin{lem}\label{fiber}
Let $E\stackrel{\pi}{\longrightarrow} M$ be a smooth fiber bundle with compact fibre.
If $M$ admits an open book decomposition, then so does $E$.
\end{lem}
\begin{proof}
If $p\colon M\longrightarrow\C$ defines an Alexander decomposition of $M$, then 
$p\circ \pi\colon E \longrightarrow\C$ defines one for $E$.
\end{proof}

Let's see what this means for surfaces. 
If we think of $S^2$ as being the unit sphere in $\R^3$, then
\begin{alignat*}{1}
p\colon S^2 &\longrightarrow \C\\
(x,y,z) &\longmapsto (x,y)
\end{alignat*}
defines an Alexander decomposition of $S^2$ with binding $B = S^0 = \{ (0,0,\pm 1)\}$.
The two-torus $T^2$
has an Alexander decomposition $p\colon T^2\longrightarrow S^1\subset \C$ which 
has empty binding $B$. It does not have one with non-empty binding because of Lemma~\ref{Euler}.

The closed oriented surfaces of genus $g\geq 2$ do not have any open book decompositions.
Again this follows from Lemma~\ref{Euler} since no collection of points can have negative
Euler characteristic. However, we can prove it in a different way, which is instructive 
for the subsequent discussion. If a surface $\Sigma_g$ had an open book decomposition,
then by Lemma~\ref{fiber} the same would be true for any fiber bundle over it. However,
since there are surface bundles over $\Sigma_g$ with non-zero signature by a classical 
construction going back to Atiyah~\cite{A} and Kodaira~\cite{Kd}, 
we would get a contradiction with Lemma~\ref{Weyl}. 

This proof extends to products in the following way.
\begin{thm}\label{main}
Let $M=\Sigma_{g_1}\times\ldots\times\Sigma_{g_k}\times N$ be a product of $k\geq 1$ closed oriented
surfaces, all of genus $g_i\geq 2$, and a closed oriented $4n$-manifold $N$ of non-zero signature. 
Then $M$ does not admit an open book decomposition.

Furthermore, for any $L$, the connected sum $M\sharp L$ does not admit an open book decomposition
either.
\end{thm}
\begin{proof}
Over each $\Sigma_{g_i}$ there exists a surface bundle 
$$
F_i\longrightarrow X_i\longrightarrow \Sigma_{g_i}
$$ 
with non-zero signature; see for example Bryan--Donagi~\cite{BD}, and the literature quoted there. 
Then 
$$
X=X_1\times\ldots\times X_k\times N
$$
is a $4(k+n)$-dimensional manifold of non-zero signature, which is a smooth fiber bundle over $M$.
The combination of Lemmata~\ref{Weyl} and~\ref{fiber} gives the conclusion for $M$.

To get the conclusion for $M\sharp L$ we fibre sum $X$ 
over $M$ to a trivial $F_1\times\ldots \times F_k$-bundle over $L$ and use Novikov additivity 
to get a bundle over $M\sharp L$ whose total space has non-zero signature.
\end{proof}

Since every surface is a double, the same is true for any product containing a surface factor.
Therefore we obtain:
\begin{cor}
In every even dimension there are doubles without open book decompositions.
\end{cor}
In dimensions $\geq 6$ this contradicts~\cite[Cor.~30.13]{Ran}. 
\begin{ex}\label{P2}
If $g(\Sigma )\geq 2$, then by Theorem~\ref{main} the $6$-manifold $\Sigma\times\C P^2$ does not have an open book decomposition,
contrary to what is claimed in~\cite[Rem.~30.14]{Ran}.
\end{ex}
As mentioned in the introduction, the proof of Theorem~\ref{main} shows that an invariant $i(M)$ whose 
vanishing is equivalent to the existence of an open book decomposition has to take into 
account the signatures of all possible fibre bundles over $M$. For example, the manifold
$\Sigma\times\C P^2$ itself has zero signature, but carries a bundle with
a total space of non-zero signature.

\section{Dimension four}

In dimension $4$, neither Quinn~\cite{Quinn} nor Ranicki~\cite{Ran} claimed to have a complete 
invariant whose vanishing is equivalent to the existence of an open book decomposition.
However, the vanishing of the Quinn invariant $i(M)$ is certainly necessary, but not known to
be sufficient.

From Theorem~\ref{main} we know that a product of two high-genus surfaces does not have an 
open book decomposition, although it has vanishing signature, and is in fact a double.
We can generate more examples using the following result.

\begin{prop}\label{p:dom}
Let $M$ be a closed oriented four-manifold admitting a map of non-zero degree to 
$\Sigma_{g_1}\times\Sigma_{g_2}$ with $g_i\geq 2$. Then $M$ carries a four-manifold
bundle $X\longrightarrow M$ with total space of non-zero signature.
\end{prop}\label{dom}
An $M$ as in the Proposition is not necessarily a double. We can replace any such $M$
by $M\sharp\overline{M}$, which is the double of $M\setminus D^4$, to generate more examples,
since if $M$ dominates $\Sigma_{g_1}\times\Sigma_{g_2}$, so does $M\sharp\overline{M}$.
\begin{proof}
We may assume that the signature of $M$ vanishes, for otherwise the claim is trivial. 
For dimension reasons this means that the Pontryagin classes of $TM$ are trivial.

For any smooth fibre bundle $\pi\colon X\longrightarrow M$ the signature of the total space is given
by Atiyah's formula~\cite{A}
\begin{equation}\label{AS}
\sigma (X) = \langle   ch (Sign(\pi ))\cdot \mathcal{\tilde{L}}(M) , [M]   \rangle \ .
\end{equation}
This is an application of the index theorem for families of signature 
operators. In the formula $\mathcal{\tilde{L}}(M)$ is a power series in the Pontryagin classes of $TM$, and $Sign(\pi ) = W_+ - W_-$
is the virtual vector bundle over the base consisting of the eigenbundles of the fibrewise Hodge star on the middle-dimensional real cohomology 
groups of the fibres of $\pi$.

We first take a bundle of non-zero signature $\pi\colon X_1\times X_2\longrightarrow\Sigma_{g_1}\times\Sigma_{g_2}$ 
as in the proof of Theorem~\ref{main} and apply this formula. Since the Pontryagin classes of the base are trivial, we 
conclude that $\mathcal{\tilde{L}}(\Sigma_{g_1}\times\Sigma_{g_2})$ vanishes in positive degrees
and the evaluation of $ch (Sign(\pi ))$ on the base gives the signature and hence  is non-zero.

Next we pull back this bundle under a map $\varphi\colon M\longrightarrow \Sigma_{g_1}\times\Sigma_{g_2}$ of 
non-zero degree. Then we have $ch (Sign (\varphi^*(\pi))) = \varphi^*(ch (Sign(\pi )))$, and since $\varphi$ has non-zero degree this
evaluates non-trivially on the base. Since the Pontryagin classes of $M$ are trivial, $\mathcal{\tilde{L}}(M)$ also vanishes in
positive degrees,
and the non-zero evaluation of the Chern character on $M$ gives the signature of the pulled back bundle.
\end{proof}
The slightly roundabout argument is needed because we do not assume that the dominant map $\varphi$ is  tangential.

\begin{ex}\label{Gaifullin}
By the work of Gaifullin~\cite{Gai} there are real hyperbolic four-manifolds that dominate products
$\Sigma_{g_1}\times\Sigma_{g_2}$ with $g_i\geq 2$ with non-zero degree. By Proposition~\ref{p:dom} these
hyperbolic four-manifolds carry bundles with non-zero signature, and so do not have an open book decompositions, 
although the signatures of the hyperbolic manifolds themselves are zero.
\end{ex}

\begin{ex}\label{Stover}
By the work of Stover~\cite{Sto}, 
$\Sigma_{g_1}\times\Sigma_{g_2}$ has a complex structure such that a suitable blowup of this product is biholomorphic to 
the analytic space underlying an orbifold ball quotient $\C H^2/\Gamma$. 
Since the lattice $\Gamma$ is virtually torsion-free
by Selberg's Lemma, one can take a finite (orbifold) covering which is a manifold ball
quotient $M=\C H^2/\Delta$, rather than an orbifold. The composition of the covering map with the blowdown
from the orbifold to the product of curves is a holomorphic map $f\colon M\longrightarrow \Sigma_{g_1}\times\Sigma_{g_2}$
of positive degree. 

Now $M$, being a ball quotient manifold, has positive signature, equal to $k>0$ say, and so does not have an open book
decomposition. However, the map $f$ shows that if we take a $k$-fold blowup $M_k = M\sharp k\overline{\mathbb{C} P^2}$
or we take $M\sharp\overline{M} = D(M\setminus D^4)$, then these manifolds of zero signature have maps of non-zero degree
to $\Sigma_{g_1}\times\Sigma_{g_2}$, and so carry bundles of non-zero signature. Therefore they do not 
have open book decompositions.
\end{ex}
These examples exhibit many different geometric properties. The real hyperbolic examples are not
homotopy equivalent to complex surfaces~\cite{Kot92}, whereas the blowups $M_k$ of a ball quotient 
are complex K\"ahler and have non-vanishing Seiberg--Witten invariants. The doubles $M\sharp\overline{M}$ 
are not complex and have trivial Seiberg--Witten invariants.
Since for all these manifolds the absence of open book decompositions is detected by the existence of 
bundles of non-zero signature, they should all have non-vanishing Quinn invariant. In fact, the argument 
above saying that manifolds dominating a product of surfaces cannot have open book decompositions should 
be interpreted as a statement about the functorial nature of the Quinn invariant.

\section{Relation with the simplicial volume}

The simplicial volume $\vert\vert M\vert\vert$ of a closed oriented manifold was defined by 
Gromov~\cite{Gro} as the $\ell^1$-norm of the fundamental class: $\vert\vert M\vert\vert = \vert\vert [M]\vert\vert_{\ell^1}$.
It is multiplicative in finite unramified covers, and satisfies 
\begin{equation}\label{eq:mult}
\vert \vert N\vert\vert \geq \vert d\vert \cdot \vert\vert M\vert\vert
\end{equation}
whenever there is a degree $d$ map $\varphi\colon N\longrightarrow M$. The argument from the proof
of Lemma~\ref{Weyl} immediately implies:
\begin{prop}
If $M$ admits an open book decomposition with monodromy $f$ of finite order, then $\vert\vert M\vert\vert =0$.
\end{prop}
\begin{proof}
Every open book with monodromy of finite order admits a map of non-zero degree, in fact, a branched covering map, from an open book 
with trivial monodromy. So it is enough to prove vanishing in the case that $f$ is the identity. In this case
$M$ maps to itself with arbitrarily large degree, so~\eqref{eq:mult} implies the conclusion.
\end{proof}

\begin{ex}
It is well known that in dimension $3$ an open book with trivial monodromy is a connected sum of 
copies of $S^1\times S^2$, so the vanishing of the simplicial volume is clear.
\end{ex}

Recently Kastenholz proved that in dimension four the assumption about the monodromy is 
not needed.
\begin{thm}[Kastenholz~\cite{Kas}]\label{t:K}
If a closed oriented four-manifold $M$ admits an open book decomposition, then $\vert\vert M\vert\vert =0$.
\end{thm}
This result is special to dimension $4$ (and $2$). In dimension $3$ all manifolds, including, say, hyperbolic
ones, have open book decompositions by the classical work of Alexander. Taking a product of such a manifold with any hyperbolic $(n-3)$-manifold
produces open books with positive simplicial volume in all dimensions $n\geq 5$ (use Lemma~\ref{fiber} for a product bundle).

Since compact real or complex hyperbolic manifolds have positive simplicial volume, the
Examples~\ref{Gaifullin} and \ref{Stover} can be generalised. By Theorem~\ref{t:K}
one can dispense with the condition that there be a dominating map to a product of surfaces,
and consider arbitrary hyperbolic manifolds. To discuss blowups and doubles of $M\setminus D^4$
one can use the additivity of the simplicial volume in connected sums~\cite{Gro}.

In a similar vein, instead of considering 
 the product $\Sigma_{g_1}\times\Sigma_{g_2}$ as in~\cite{Kas},
 one can consider any surface bundle $M$ over a surface, with base and fibre both
of genus $\geq 2$. Then $M$ has non-zero simplicial volume by~\cite{HK},
so Theorem~\ref{t:K} applies to rule out open book decompositions. 
Of course it can happen that $N$ has non-zero signature, in which case 
there is nothing to prove because of  Lemma~\ref{Weyl}.

We now prove that all examples of four-manifolds with zero signature carrying bundles of non-zero signature 
can also be detected by Theorem~\ref{t:K}.
\begin{prop}\label{p:simpl}
Let $M$ be a closed oriented four-manifold of zero signature carrying a bundle $\pi\colon X\longrightarrow M$
with total space $X$ of non-zero signature. Then the simplicial volume of $M$ is positive: $\vert\vert M\vert\vert > 0$.
\end{prop}
In the special situation of Proposition~\ref{p:dom} this follows from~\eqref{eq:mult}. 
The point here is to prove positivity of the simplicial volume in greater
generality.
\begin{proof}
Again we use Atiyah's formula~\eqref{AS} for the signature of $X$. Since $M$ has zero signature, its Pontryagin
classes are zero, and so~\eqref{AS} simplifies since the term $\mathcal{\tilde{L}}(M)$ vanishes in positive degrees.
The assumption $\sigma (X)\neq 0$ implies that the degree $4$ part $ch_2(Sign(\pi ))$ of the Chern character of 
the virtual bundle $sign (\pi ) = W_+-W_-$ evaluates non-trivially on the fundamental class $[M]$. This signature 
bundle is a pullback from $BO(p,q)^{\delta}$ under the map corresponding to the representation
$$
\rho\colon\pi_1(M)\longrightarrow O(p,q) \ ,
$$
given by the middle-dimensional cohomology of the fibre of $\pi$. (The signature of this fibre is $p-q$.) 
Therefore, the characteristic class $ch_2(Sign(\pi ))$
is a bounded cohomology class, see Gromov~\cite[Section~1.3]{Gro}. Since its evaluation on $[M]$ is non-zero,
we conclude 
$$
\vert\vert M\vert\vert = \vert\vert [M]\vert\vert_{\ell^1} > 0
$$ 
by the duality of the sup norm on cohomology and the $\ell^1$-norm on homology.
\end{proof}

\begin{cor}\label{cc}
Let $M$ be a closed oriented four-manifold with a smooth fibre bundle $\pi\colon X\longrightarrow M$ with fibre $F$ 
such that $\sigma (X) \neq \sigma (F)\cdot \sigma(M)$.
Then $\vert\vert M\vert\vert > 0$.
\end{cor}
\begin{proof}
We can fibre sum $X$ to the product bundle $\overline{M}\times F$ to get a bundle over $M\sharp\overline{M}$
with signature
$$
\sigma (X)+\sigma (\overline{M})\cdot\sigma (F) = \sigma (X)-\sigma (M)\cdot\sigma (F) \neq 0 \ .
$$
Since $M\sharp\overline{M}$ has zero signature, Proposition~\ref{p:simpl} implies 
$$
 \vert\vert M\vert\vert  =\frac{1}{2}\cdot\vert\vert M\sharp\overline{M}\vert\vert > 0 \ .
$$
\end{proof}
In Example~\ref{P2} we mentioned the six-manifold $\Sigma\times\C P^2$ with  $g(\Sigma)\geq 2$.
This carries a bundle of non-zero signature, but has vanishing simplicial volume. Together with analogous examples
in higher dimensions, this shows that Proposition~\ref{p:simpl} has no analog in any higher dimension.
Like Theorem~\ref{t:K}, it is special to dimension $4$.
Curiously, in dimension $2$ the analogous statement is true: a closed oriented surface has non-zero simplicial 
volume if and only if it carries a bundle of non-zero signature, and this property characterises the absence of 
open book decompositions.
One might speculate that the nice situation for surfaces continues to hold for four-manifolds, so that Corollary~\ref{cc}
admits a converse:
\begin{conj}
Let $M$ be a closed oriented four-manifold with $\vert\vert M\vert\vert > 0$. Then there is a smooth
fibre bundle $\pi\colon X\longrightarrow M$ with fibre $F$ such that $\sigma (X) \neq \sigma (F)\cdot \sigma(M)$.
\end{conj}
In other words, for a four-manifold to have non-zero simplicial volume is equivalent to there being a fibre
bundle with non-multiplicative 
signature over it.
Maybe the following class of examples can serve as an interesting testing ground for the conjecture.
\begin{prop}\label{RTKW}
There are integral homology four-spheres which are doubles and have positive simplicial volume.
\end{prop}
\begin{proof}
Ratcliffe and Tschantz~\cite{RT} constructed integral homology four-spheres $M$ carrying metrics of 
non-positive sectional curvature. By the recent result of Kim and Wan~\cite{KimWan}, these manifolds have
positive simplicial volume. For every such $M$ the connected sum $M\sharp\overline{M}$ is a double
and is again an integral homology sphere. By the additivity of the simplicial volume in connected sums~\cite{Gro} 
we have $\vert \vert M\sharp\overline{M}\vert\vert >0$.
\end{proof}

\section{Final remarks}

\subsection{Application to Engel structures}

Recall that an Engel structure on a $4$-manifold $M$ is a maximally non-integrable rank $2$ subbundle 
$\mathcal{D}\subset TM$. This means that 
$\mathcal{E}=[ \mathcal{D} , \mathcal{D}]$ has rank $3$ and $[\mathcal{D},\mathcal{E}]=TM$.
Under suitable orientability assumptions, the existence of an Engel structure implies that $M$ is 
parallelizable. It was proved by Vogel~\cite{Vog} that conversely every closed parallelizable $4$-manifold 
does indeed admit an Engel structure. Colin, Presas and Vogel~\cite{CPV} introduced the notion
of a supporting open book decomposition for an Engel structure, and asked whether every Engel
structure (up to homotopy) is supported by an open book decomposition. 

The results of this paper provide many examples of parallelizable four-manifolds which do not 
admit any open book decomposition. Since they all have Engel structures~\cite{Vog}, this gives
a negative answer to the question of Colin--Presas--Vogel~\cite{CPV}. 

This negative answer was also noticed independently by Lawande and Saha~\cite{LS}, 
based on the work of Kastenholz~\cite{Kas}, who had pointed out that for $g_i\geq 2$ 
the product $\Sigma_{g_1}\times\Sigma_{g_2}$ has positive simplicial volume and therefore 
has no open book decomposition. Connect summing 
this product with the correct number of copies of $S^1\times S^3$ one gets a parallelizable 
four-manifold, and Lawande--Saha~\cite{LS} inferred that this has no open book decomposition
by the additivity of the simplicial volume combined with Theorem~\ref{t:K}. For this particular 
example one can use Theorem~\ref{main} and Proposition~\ref{dom} instead of Theorem~\ref{t:K},
so no discussion of the simplicial volume is needed. 

In the same way one can treat the hyperbolic manifolds from Example~\ref{Gaifullin}. They
have zero signature and positive even Euler characteristic~\cite{KotMZ}, and, after passing to a 
suitable finite covering, they are spin~\cite{DS}. Therefore the connected sum with the correct 
number of copies of $S^1\times S^3$ is parallelizable, and therefore~\cite{Vog} Engel. However,
there is no open book decomposition by Theorem~\ref{main} and Proposition~\ref{dom}.

If one uses the simplicial volume, then many other examples can be constructed. In particular,
Proposition~\ref{RTKW} together with Theorem~\ref{t:K} imply:
\begin{cor}
There are integral homology spheres $M$ with the property that for any closed parallelizable 
$4$-manifold $P$ the connected sum $P\sharp M$  has no open book decomposition.
\end{cor}
Note that $P\sharp M$ is again parallelizable, and is $\Z$-homology equivalent to $P$.

\subsection{Relation with the SK-groups}

The fact that products of surfaces do not have open book decompositions, proved here as part of Theorem~\ref{main},
was clearly known (to some people) in the 1970s. With hindsight, it appears in Neumann's paper~\cite{Neu}.
The argument given there, in the context of the cutting and pasting groups, or SK-groups~\cite{KKNO},
is the same as ours. It is interesting to note that Neumann's paper preceded Quinn's~\cite{Quinn} by four years,
and so Quinn could or should have known that there are doubles without open book decompositions.
Ranicki's book~\cite{Ran} in fact refers to both~\cite{Neu} and~\cite{KKNO}. I know from the published 
record, and from personal conversations, that in later years, long after the publication of~\cite{Ran},
Ranicki was very much aware of the phenomenon of non-multiplicativity of the signature in fibre bundles,
which is the obstruction to having open book decompositions on products of surfaces. Why he never made 
the connection to the discussion in his book remains a mystery.

In~\cite{Ran}, Ranicki introduced an invariant
$$
\sigma^*(M)\in LAsy^n(\Z[\pi_1(M)]) \ ,
$$
which he called the asymmetric signature. It seems he never really said explicitly that this was equivalent 
to the Quinn invariant $i(M)$, but it was clearly meant to be. In Chapters~29 and~30 of~\cite{Ran},
Ranicki concluded that in high dimensions the vanishing of the asymmetric signature is, on the one hand,
equivalent to the existence of an open book decomposition, and, on the other, equivalent to a decomposition
as a twisted double. This then led him to the erroneous conclusion that every twisted double is an open book.
The list of errata to~\cite{Ran} posted at~\cite{RanErr} makes no mention of this problem.

\bigskip


\bigskip

\bibliographystyle{amsplain}

\end{document}